\documentclass[12pt,twoside,a4paper]{article}

\usepackage[T1]{fontenc}
\usepackage[utf8]{inputenc}
\usepackage[english]{babel}
\usepackage{amsmath,amsfonts,amsthm,amssymb,latexsym}
\usepackage{geometry}
\renewcommand{\phi}{\varphi}
\renewcommand{\epsilon}{\varepsilon}

\theoremstyle{plain}
\newtheorem{theorem}{Theorem}

\theoremstyle{definition}
\newtheorem{definition}{Definition}

\theoremstyle{remark}

\newtheorem*{theorem*}{Theorem}
\newtheorem*{proposition*}{Proposition}
\newtheorem*{lemma*}{Lemma}
\newtheorem*{corollary*}{Corollary}
\newtheorem*{remark*}{Remark}
\newtheorem*{definition*}{Definition}
\newtheorem*{example*}{Example}
\newtheorem*{examples*}{Examples}

\makeatletter

\newcommand{\shorttitle}[1]{}
\newcommand{\pyear}[1]{}
\newcommand{\volume}[1]{}
\newcommand{\issue}[1]{}
\newcommand{\pageno}[1]{\setcounter{page}{#1}}
\newcommand{\UDC}[1]{}
\newcommand{\received}[1]{}

\gdef\@institute{}
\newcommand{\institute}[1]{\gdef\@institute{#1}}

\gdef\@authoremail{}
\newcommand{\email}[1]{\gdef\@authoremail{#1}}

\gdef\@subjclass{}
\newcommand{\subjclass}[1]{\gdef\@subjclass{#1}}

\long\gdef\@enabstracttext{}
\gdef\@enkeywords{}
\newcommand{\enabstract}[4]{%
  \long\gdef\@enabstracttext{#3}%
  \gdef\@enkeywords{#4}%
}
\newcommand{\printabstract}{%
  \ifx\@enabstracttext\empty\else
    \begin{abstract}
      \@enabstracttext
      \par\medskip
      \noindent\textit{Key words and phrases:} \@enkeywords.
    \end{abstract}
  \fi
}

\let\orig@maketitle\maketitle
\renewcommand{\maketitle}{%
  \orig@maketitle
  \ifx\@institute\empty\else
    \begin{center}
      \small\textit{\@institute}%
      \ifx\@authoremail\empty\else
        \\\small e-mail: \textit{\@authoremail}%
      \fi
    \end{center}
  \fi
  \ifx\@subjclass\empty\else
    \begin{center}
      \small 2010 \textit{Mathematics Subject Classification:} \@subjclass.
    \end{center}
  \fi
}

\makeatother

\author{Mokliachuk Oleksandr}
\title{Estimation of reliability and accuracy of models of $\varphi$-sub-Gaussian process in $L_p(T)$ using generating functions of polynomial expansions}

\shorttitle{Estimation of reliability and accuracy of models of $\varphi$-sub-Gaussian process in $L_p(T)$}

\enabstract{Mokliachuk Oleksandr}%
{Estimation of reliability and accuracy of models of $\varphi$-sub-Gaussian process in $L_p(T)$ using generating functions of polynomial expansions}%
{Stochastic processes are often represented through orthonormal series expansions, a framework originating in the classical works of Lo\`eve and Karhunen and widely used for simulation and numerical approximation. While truncation error in such expansions has been extensively studied, practical models frequently involve an additional source of error arising from the approximation of coefficient functions when closed-form expressions are unavailable. The combined effect of these two errors remains insufficiently addressed in the literature. Building on the author's earlier work on reliability and accuracy estimates for $\varphi$-sub-Gaussian processes, this paper extends the methodology to orthonormal polynomial systems that do not possess normalized generating functions in analytical form, including the Legendre, generalized Laguerre, and Gegenbauer families. New bounds are derived for models in $L_p(T)$ space that simultaneously account for truncation and coefficient approximation. The resulting criteria provide practical guidance for selecting the number of series terms required to achieve prescribed levels of reliability and accuracy across a broader class of polynomial-based stochastic process models.
}
{Models of stochastic processes, quasi-Banach spaces, Karhunen-Lo{\`e}ve model, reliability and accuracy, sub-Gaussian random processes,
 Legendre polynomials, generalized Laguerre polynomials, Gegenbauer polynomials}


\UDC{519.2}

\pyear{2026}
\volume{14}
\issue{1}
\pageno{10}
\received{01.12.2025}

\institute{Sun Pharma Advanced Research Company, Ltd. -- Princeton, NJ, USA}
\email{omokliachuk@gmail.com (Mokliachuk O.M.)}

\begin{document}

\maketitle
\printabstract


\section*{Introduction}
Stochastic processes arise in a wide range of applications in applied probability, statistical modeling, signal processing, and quantitative finance. In many of these settings, one requires a tractable model that allows simulation or numerical approximation of a given process. A standard and powerful approach is to represent the process through an expansion over an orthonormal functional system, $X(t) = \sum_{k=0}^{\infty} a_k(t)\,\xi_k$, as established in the foundational works of Lo\`eve and Karhunen \cite{loeve,Karhunen1947}. Truncated versions of such expansions form the basis of numerous simulation and approximation techniques; see, for example, \cite{AsmussenGlynn2007,GhanemSpanos1991}. Classical decomposition theorems, including those in \cite{kozroztur}, justify these representations and describe how the covariance structure of a process determines the coefficient functions $a_k(t)$.

Truncating the expansion naturally introduces error, and quantifying this truncation error has been addressed in many papers, for example, in \cite{koz-models}. However, in many applications the coefficients $a_k(t)$ cannot be computed analytically and must instead be replaced by approximations $\hat a_k(t)$. This second source of error - arising from coefficient approximation - is generally not treated in the existing literature and requires a separate analysis. In practical modeling scenarios, both effects occur together and must be incorporated into a unified reliability and accuracy analysis.

In the author's previous publications, a general method was developed for obtaining reliability and accuracy estimates for models of $\varphi$-sub-Gaussian processes in $L_p(T)$ (see \cite{KMM2015,Chapter,Mokliachuk2012,Mokliachuk2014}) and in $C(T)$ (see \cite{Chapter,Mokliachuk2018}) under the assumption that the underlying orthonormal system admits a generating function. This framework enabled the analysis of tail behavior for truncated orthonormal expansions and the derivation of explicit upper bounds for modeling errors. Classical systems such as Hermite and Chebyshev polynomials served as illustrative examples; however, many important orthogonal polynomial families do not possess normalized generating functions in closed analytical form.

The present paper extends this methodology to a broader class of polynomial systems. In particular, we establish reliability and accuracy estimates in $L_p(T)$ for models based on Legendre, generalized Laguerre, and Gegenbauer orthonormal polynomial families, while remaining within the general framework of \cite{KMM2015,Chapter,Mokliachuk2012,Mokliachuk2014}.
These results yield practical criteria for selecting the number of series terms
required to achieve prescribed levels of reliability and accuracy.

The structure of the paper is as follows. Section 2 introduces the basic notation and recalls the decomposition theorem from \cite{kozroztur}. Section 3 reviews key properties of $\varphi$-sub-Gaussian random variables and summarizes earlier results for models in $L_p(T)$. Section 4 derives new bounds for above three orthonormal polynomial families.
The paper concludes with a discussion of applicability and future directions.

\section{Basic notations}
Consider a second-order stochastic process  $X=\{X(t),t\in \mathbf{T}\}$.
 As proved in \cite{kozroztur}, the next statement holds true.

\begin{theorem}\label{theor_krt} \rm{(On decomposition of the stochastic
process using an orthonormal basis.)}
Assume $ X(t),t\in \mathbf{T}$ is a centered second order stochastic process, $EX(t) = 0$, $E|X(t)|^2<\infty$, $t \in \mathbf{T}$.
Let $B(t, s) = EX(t)\overline{X(s)}$ be the correlation function of the process $X(t)$.
Assume $f(t, \lambda)$ is a function from the space $L_2(\Lambda, \mu)$ for all
$t \in \mathbf{T}$,
and
$\{g_k(\lambda),k=0,1,\dots\}$
 is an orthonormal basis in this space.

The correlation function $B(t, s)$ can be represented as
$$B(t,s)=\int_\Lambda f(t,\lambda)\overline{f(s,\lambda)}d\mu(\lambda)$$
if and only if the following representation holds true
\begin{equation}
\label{mainseries}
X(t)=\sum_{k=0}^\infty a_k(t)\xi_k,
\end{equation}
where $$a_k(t)=\int_\Lambda f(t,\lambda)\overline{g_k(\lambda)}d\lambda.$$
Here, $\xi_k$ are centered uncorrelated random variables: $E\xi_k=0$, $E\xi_k\overline{\xi_l}=\delta_{kl}$,
$E|\xi_k|^2=1$.
\end{theorem}

In \cite{koz-models}, \cite{Chapter}, the model of stochastic process is considered as the sum of the first $N$ elements of the decomposition (\ref{mainseries}).
However, for a significant number of processes, except for a very limited class, it is difficult or impossible to find $a_k(t)$ explicitly.
Therefore, in this paper we will assume that it is possible to find some approximations $\hat{a}_k(t)$ of the functions $a_k(t)$.
Since approximation of $a_k(t)$ will contribute to the overall modeling error, it is important to include this effect in the overall reliability and accuracy estimate of the model. Thus, in this paper we will use the following notation.

\begin{definition}
For a stochastic process $X(t), t\in\mathbf{T},$ that allows representation (\ref{mainseries}),  we call the process
\begin{equation}
\label{model}
X_N(t)=\sum_{k=0}^N \xi_k \hat{a}_k(t),\,\, t\in\mathbf{T},
\end{equation}
a model of the process $X(t), t\in\mathbf{T}$.
Here, $\hat{a}_k(t)$ are approximations of functions ${a}_k(t)$ in  representation (\ref{mainseries}), random variables $\xi_k$ satisfy conditions of Theorem \ref{theor_krt}: $E\xi_k=0$, $E\xi_k\overline{\xi_l}=\delta_{kl}$,
$E|\xi_k|^2=1$.
\end{definition}

It will be useful to introduce the following notation: $$\Delta_N(t)=X(t)-X_N(t).$$

Following the statement of the Theorem \ref{theor_krt}, we may consider a system of orthonormal functions $\{g_k(\lambda)\}$ as an orthonormal basis that can be used to represent a stochastic process using (\ref{mainseries}). A good candidate for this role may be an orthonormal polynomial system, since these possess some useful properties that can be used for estimations of reliability and accuracy of models of stochastic processes.

One useful property of systems of orthogonal polynomials is the existence of a generating function.

\begin{definition}
For a system of orthogonal polynomials $g_k(t)$, a generating function is a sum of a formal series
$$GF(t,w)=\sum_{k=0}^\infty g_k(t)w^k.$$
It is assumed that $w$ is chosen so that the series converges.
\end{definition}

For a lot of classic orthogonal polynomial systems, such as Cauchy, Hermite, Legendre, and many more, generating functions are well known. However, for use as a decomposition basis for the process, we require not just an orthogonal system, but an orthonormal one. The derivation of the generation function for an orthonormal version of classic orthogonal polynomials may not be straightforward, and additional conditions may be required to be introduced to the final results.

\section{Model estimations in $L_p[0,T]$}

Let  $\left(\Omega,{\cal F}, P \right)$ be a standard probability space, let $ L_2(\Omega)$ be the space of centered random variables with bounded second moments, $E\xi^2<\infty$. Assume that $\{\Lambda, {\cal U}, \mu\}$ is a measurable space with a $\sigma$-finite measure $\mu$. Let $L_p(\Lambda,\mu)$ be a Banach space of integrable functions in the power $p$ with measure $\mu$.

\begin{definition}
\cite{Buldygin2000,koz-models} A random variable $\xi$ is called sub-Gaussian if there exist $a\geq 0$ such that for all $\lambda\in R$ the inequality
$$E\exp\{\lambda\xi\}\leq\exp\left\{\frac{a^2\lambda^2}{2}\right\}$$
holds true.

\noindent The function
$$\tau(\xi)=\inf\left\{a\geq0: E \exp\{\lambda\xi\}\leq \exp\left\{\frac{a^2\lambda^2}{2}\right\},\lambda\in R\right\}$$
is a sub-Gaussian norm of a random variable $\xi$.
\end{definition}

\begin{definition}
A continuous even convex function $\varphi=\{\varphi(x),x\in R\}$ is called an Orlicz $N$-function, if $\varphi(0)=0$ and $\varphi(x)>0$ for $x\neq0$ and

$$\lim_{x\to0}\frac{\varphi(x)}{x}=0 \mbox{  and  } \lim_{x\to\infty}\frac{\varphi(x)}{x}=\infty.$$
\end{definition}

\begin{definition}
Let $\varphi=\{\varphi(x),x\in R\}$ be an Orlicz $N$-function, and let $$\lim_{x\to 0} \inf_x \frac{\varphi(x)}{x^2}=c>0.$$
A random variable $\xi$ belongs to the space $Sub_\varphi(\Omega)$, if $E\xi=0$, $E\exp\{\lambda\xi\}$ exists for every $\lambda\in R$,
and there exists a constant $a>0$ such that for all $\lambda\in R$ the inequality
$$E\exp\{\lambda\xi\}\leq \exp\{\varphi(\lambda a)\}$$ holds true.
\end{definition}

In \cite{Kozachenko1998} it is proved that the space $Sub_\varphi(\Omega)$ is a Banach space with respect to the norm

$$\tau_\varphi(\xi)=\sup_{\lambda>0} \frac{\varphi^{(-1)} (\ln\  E \exp\{\lambda\xi\})}{\lambda}$$

\begin{definition} A stochastic process $X=\{X(t),t\in T\}$ is called $\varphi$-sub-Gaussian, if random variables $X(t)$ for each $t\in T$ are $\varphi$-sub-Gaussian.
\end{definition}

\begin{definition}
Let $X=\{X(t),t\in [0,T]\}$ be a stochastic process belonging to the space $Sub_\varphi(\Omega)$.
We will call stochastic process $X_N=\{X_N(t),t\in [0,T]\}$ from $Sub_\varphi(\Omega)$ a model that approximates $X$ with giver reliability $1-\alpha$ and accuracy $\delta$ in the space $L_p[0,T]$ if
$$P\left\{\left(\int_0^T( X(t)-X_N(t))^pdt\right)^{1/p}>\delta\right\}\leq\alpha.$$
\end{definition}

The following two theorems were proved in \cite{Chapter}.

\begin{theorem}
\label{lp_l2}
 Let a stochastic process $X=\{X(t),t\in [0,T]\}$ belong to the space ${Sub}_\varphi (\Omega)$ and let
  $$\varphi(t)=\frac{t^\gamma}{\gamma},\  1<\gamma\leq 2.$$
  Assume that $$C_N=\int_0^T (\tau_\varphi(X(t)-X_N(t)))^pd\mu(t)<\infty.$$

\noindent The model $X_N(t)$ approximates the stochastic process $X(t)$ with given reliability $1-\alpha$ and accuracy $\delta$ in the space
$L_p[0,T]$,
if

$$\left\{\begin{array}{c}C_N\leq \delta/(\beta \ln \frac{2}{\alpha})^{p/\beta},\\  C_N<\delta /p^{p\left(1-1/\gamma\right)},\end{array}\right.$$
where $\beta$ if a number that fits the condition $\frac{1}{\beta}+\frac{1}{\gamma}=1$.

\end{theorem}

\begin{theorem}
\label{lp_g2}
  Let a stochastic process $X=\{X(t),t\in [0,T]\}$ belong to the space $\text{S u b}_\varphi (\Omega)$ and
  let
  $$\varphi(t)=\left\{\begin{array}{c}\frac{t^2}{\gamma}, t<1,\\   \frac{t^\gamma}{\gamma},t\geq1,\end{array}\right.$$
where $\gamma>2$. Assume that $$C_N=\int_0^T (\tau_\varphi(X(t)-X_N(t)))^pd\mu(t)<\infty.$$

\noindent The model $X_N(t)$ approximates the stochastic process $X(t)$ with given reliability $1-\alpha$ and accuracy $\delta$ in the space
$L_p[0,T]$,
 if
$$\left\{\begin{array}{c}C_N\leq \delta/(\beta \ln \frac{2}{\alpha})^{p/\beta}, \\  C_N<\delta/p^{p(1-1/\gamma)},\end{array}\right.$$
where $\frac{1}{\beta}+\frac{1}{\gamma}=1$.

\end{theorem}

It follows from  Theorem \ref{lp_l2} and  Theorem \ref{lp_g2}, that in order to find the required number of series elements $N$ we  have to estimate the constant
$$C_N = \int_0^T \left(\tau_\varphi\left(X(t)-X_N(t)\right)\right)^pd\mu(t)$$
that depends on the difference between the model of the process and the process itself.

The key element of this expression is $\tau_\varphi\left(X(t)-X_N(t)\right)$, which can be estimated as following:

$$\tau_\varphi(\Delta_N(t)) = \tau_\varphi(X(t)-X_N(t))=\tau_\varphi \left(\sum_{k=0}^\infty \xi_k a_k(t) - \sum_{k=0}^N \xi_k \hat{a}_k(t)\right) = $$

$$=\tau_\varphi \left(\sum_{k=0}^N \xi_k \delta_k(t) + \sum_{k=N+1}^\infty \xi_k {a}_k(t)\right) \leq \sum_{k=0}^N \tau_\varphi(\xi_k) \delta_k(t) +\sum_{k=N+1}^\infty \tau_\varphi(\xi_k) {a}_k(t) = $$
$$=\sum_{k=0}^\infty \tau_\varphi(\xi_k) a_k(t) - \sum_{k=0}^N \tau_\varphi(\xi_k) \hat{a}_k(t).$$

This inequality can be used in the results that follow.

In the next subsections we will consider several orthonormal polynomial systems that can serve as a basis for the process decomposition, namely, Legendre, generalized Laguerre, and Gegenbauer polynomials.

\subsection{Legendre polynomials}

The Legendre polynomials $P_k(t)$ are defined for $t\in[-1;1]$ as solutions to the Legendre differential equation

$$(1-t^2)y''-2yy'+k(k+1)y=0,\ y=y(t).$$

\noindent The Rodriguez formula for Legendre polynomials if the following:

$$P_k(t)=\frac{1}{2^k k!} \frac{d^k}{dt^k}(t^2-1)^k.$$

\noindent  The generation function for the Legendre polynomials is the following:

$$GF_{Legendre}(t,w)=\sum_{k=0}^\infty P_k(t)w^k = \frac{1}{\sqrt{1-2 t w + w^2}}.$$

\noindent  These polynomials form the orthogonal system:

$$\int_{-1}^1 P_n(t)P_m(t)dt = \frac{2}{2n+1} \delta_{nm},$$
where $\delta_{nm}$ is the Kronecker delta.

\noindent Therefore, this polynomial system can generate the orthonormal basis $\{\hat{P}_k(t),\,t\in[-1;1]\}$, where
\begin{equation}
\label{P-legendre}
\hat{P}_k(t)=\sqrt{\frac{2k+1}{2}}P_k(t).
\end{equation}
Unlike the polynomial systems considered in \cite{Chapter}, the generating function of the system $\{\hat{P}_k(t)\}$ is hard to derive. It may even not be possible to express this generating function in elementary functions. Therefore, to compensate for elements that obstruct calculation of this function, we will introduce the upper bound on the sub-Gaussian norm in the following form:

\begin{equation}
\label{tau-legendre}
\tau_\varphi(\xi_k)\leq \tau_{Legendre}(w,k)=\sqrt{\frac{2}{2k+1}}\tau w^k,
\end{equation}
$\tau$ is a constant.
From this moment and further in the paper, we will also assume that
$\tau_\varphi (a_k(u))\leq \tau_\varphi (\hat{a_k}(u))$, $\forall u\in [-1,1]$.

To estimate $C_N$, we will apply the following estimates:
$$\tau_\varphi(\Delta_N(t))\leq  \sum_{k=0}^\infty \tau_\varphi(\xi_k) a_k(t) - \sum_{k=0}^N \tau_\varphi(\xi_k) \hat{a}_k(t)  $$
$$=\sum_{k=0}^\infty \tau_\varphi(\xi_k) \int_{-1}^1 f(t,\lambda) \hat{P}_k(\lambda) d\lambda - \sum_{k=0}^N \tau_\varphi(\xi_k) \hat{a}_k(t) $$
$$= \int_{-1}^1 \left(
 {\sum_{k=0}^{\infty}}\tau_\varphi(\xi_k)f(t,\lambda) \hat{P}_k(\lambda) \right)d\lambda- \sum_{k=0}^N \tau_\varphi(\xi_k) \hat{a}_k(t) $$
$$\leq \int_{-1}^1 \left( {\sum_{k=0}^{\infty}}\tau_\varphi(\xi_k)f(t,\lambda) \hat{P}_k(\lambda) \right)d\lambda- \sum_{k=0}^N \tau_\varphi(\xi_k) \hat{a}_k(t) $$
$$\leq \left( \int_{-1}^1 |f(t,\lambda)|^2d\lambda\right)^{1/2}\left(\int_{-1}^1 \left( {\sum_{k=0}^{\infty} }\tau_\varphi(\xi_k) \hat{P}_k(\lambda)\right)^2 d\lambda\right)^{1/2}- \sum_{k=0}^N \tau_\varphi(\xi_k) \hat{a}_k(t). $$

Applying inequality (\ref{tau-legendre}), we will have
$$\left( \int_{-1}^1 |f(t,\lambda)|^2d\lambda\right)^{1/2}\left(\int_{-1}^1 \left( {\sum_{k=0}^{\infty}} \tau_\varphi(\xi_k) \hat{P}_k(\lambda)\right)^2 d\lambda\right)^{1/2}- \sum_{k=0}^N \tau_\varphi(\xi_k) \hat{a}_k(t)  $$
$$\leq \left( \int_{-1}^1 |f(t,\lambda)|^2d\lambda\right)^{1/2}\left(\int_{-1}^1 \left( {\sum_{k=0}^{\infty}} \sqrt{\frac{2}{2k+1}}\tau w^k \hat{P}_k(\lambda)\right)^2 d\lambda\right)^{1/2}- \sum_{k=0}^N \sqrt{\frac{2}{2k+1}}\tau w^k \hat{a}_k(t) $$
$$= \left( \int_{-1}^1 |f(t,\lambda)|^2d\lambda\right)^{1/2}\left(\int_{-1}^1 \left( {\sum_{k=0}^{\infty}} \sqrt{\frac{2}{2k+1}}\tau w^k \sqrt{\frac{2k+1}{2}}P_k(\lambda)\right)^2 d\lambda\right)^{1/2}
$$
$$- \sum_{k=0}^N \sqrt{\frac{2}{2k+1}}\tau w^k \hat{a}_k(t) $$

Canceling the square roots and substituting the first series with the corresponding generating function, we obtain
$$ {\tau} \left( \int_{-1}^1 |f(t,\lambda)|^2d\lambda\right)^{1/2}\left(\int_{-1}^1 \left( {\sum_{k=0}^{\infty}} w^k P_k(\lambda)\right)^2 d\lambda\right)^{1/2}- \sum_{k=0}^N \sqrt{\frac{2}{2k+1}}\tau w^k \hat{a}_k(t)
$$
$$= \tau \left( \int_{-1}^1 |f(t,\lambda)|^2d\lambda\right)^{1/2}\left(\int_{-1}^1 \left(\frac{1}{\sqrt{1-2 \lambda w + w^2}}\right)^2 d\lambda\right)^{1/2}- \sum_{k=0}^N \sqrt{\frac{2}{2k+1}}\tau w^k \hat{a}_k(t)
$$
$$= \tau \left( \int_{-1}^1 |f(t,\lambda)|^2d\lambda\right)^{1/2}\left(\int_{-1}^1 \frac{1}{1-2 \lambda w + w^2}d\lambda\right)^{1/2}  - \sum_{k=0}^N \sqrt{\frac{2}{2k+1}}\tau w^k \hat{a}_k(t)
$$
$$= \tau \left( \int_{-1}^1 |f(t,\lambda)|^2d\lambda\right)^{1/2}\left(\left.\frac{\ln(1-2\lambda w + w^2)}{2w}\right|_{-1}^1\right)^{1/2}  - \sum_{k=0}^N \sqrt{\frac{2}{2k+1}}\tau w^k \hat{a}_k(t)
$$
$$= \tau \left( \int_{-1}^1 |f(t,\lambda)|^2d\lambda\right)^{1/2}\frac{1}{\sqrt{w}}\sqrt{\ln\left(\frac{w+1}{w-1}\right)}  - \sum_{k=0}^N \sqrt{\frac{2}{2k+1}}\tau w^k \hat{a}_k(t). $$

Therefore, we have the following estimate for the $C_N$:

$$C_N\leq {\int_0^T }\left(\tau \left( \int_{-1}^1 |f(t,\lambda)|^2d\lambda\right)^{1/2}\frac{1}{\sqrt{w}}\sqrt{\ln\left(\frac{w+1}{w-1}\right)}  - \sum_{k=0}^N \sqrt{\frac{2}{2k+1}}\tau w^k \hat{a}_k(t)\right)^p d\mu(t).$$

Let us denote the right side of the last inequality by
$C_{N,Legendre}$.

It follows from  Theorem \ref{lp_l2} and  Theorem \ref{lp_g2}, that the following statements are true.

\begin{theorem}
\label{Legendre1}
 Let a stochastic process $X=\{X(t),t\in [0,T]\}$ belong to the space ${Sub}_\varphi (\Omega)$ with the Orlicz $N$-function
  $$\varphi(t)=\frac{t^\gamma}{\gamma},\,\,1<\gamma\leq 2,$$
  and let the process $X(t)$ admits the orthogonal decomposition  (\ref{mainseries})
 based on Legendre orthonormal polynomial families (\ref{P-legendre}).
Assume that $C_{N,Legendre}<\infty$,
$\tau_\varphi (a_k(u))\leq \tau_\varphi (\hat{a}_k(u))$, $\forall u\in [0,\infty)$ and condition \eqref{tau-legendre} holds true.
The model (\ref{model}) $X_N(t)=\sum_{k=0}^N \xi_k \hat{a}_k(t)$ approximates the stochastic process $X(t)$
with given reliability $1-\alpha$ and accuracy $\delta$ in the space $L_p[0,T]$, if
$$\left\{\begin{array}{c}C_{N,Legendre}\leq \delta/(\beta \ln \frac{2}{\alpha})^{p/\beta},\\  C_{N,Legendre}<\delta /p^{p\left(1-1/\gamma\right)},\end{array}\right.$$
where $\beta$ if a number that fits the condition $\frac{1}{\beta}+\frac{1}{\gamma}=1$.
\end{theorem}

\begin{theorem}
\label{Legendre2}
 Let a stochastic process $X=\{X(t),t\in [0,T]\}$ belong to the space ${Sub}_\varphi (\Omega)$ with the Orlicz $N$-function
   $$\varphi(t)=\left\{\begin{array}{c}\frac{t^2}{\gamma}, t<1,\\  \frac{t^\gamma}{\gamma},t\geq1,\end{array}\right.,$$
 where $\gamma>2$,
  and let the process $X(t)$ admits the orthogonal decomposition  (\ref{mainseries})
 based on Legendre orthonormal polynomial families (\ref{P-legendre}).
Assume that $C_{N,Legendre}<\infty$,
$\tau_\varphi (a_k(u))\leq \tau_\varphi (\hat{a}_k(u))$, $\forall u\in [-1,1]$ and condition \eqref{tau-legendre} holds true.
 The model (\ref{model}) $X_N(t)=\sum_{k=0}^N \xi_k \hat{a}_k(t)$ approximates the stochastic process $X(t)$
with given reliability $1-\alpha$ and accuracy $\delta$ in the space $L_p[0,T]$, if
$$\left\{\begin{array}{c}C_{N,Legendre}\leq \delta/(\beta \ln \frac{2}{\alpha})^{p/\beta},\\  C_{N,Legendre}<\delta /p^{p\left(1-1/\gamma\right)},\end{array}\right.$$
where $\beta$ if a number that fits the condition $\frac{1}{\beta}+\frac{1}{\gamma}=1$.
\end{theorem}

\subsection{Generalized Laguerre polynomials}

The Laguerre polynomials $L_k(t)$ are defined for $t\in[0;\infty]$ as solutions of the Laguerre differential equation

$$ty''+(1-t)y'+ky=0, y=y(t),$$
where $k$ is a non-negative integer. The generalized Laguerre polynomials are solutions of the following equation:

$$ty''+(\alpha + 1-t)y'+ky=0,\ y=y(t).$$
Here, an arbitrary real $\alpha$ is introduced in the equation. Hence, the generalized Laguerre polynomials are denoted as $L_k^{(\alpha)}(t)$.

\noindent The Rodriguez formula for generalized Laguerre polynomials has the form:

$$L_k^{(\alpha)}(t)=\frac{t^{-\alpha}}{k!}\frac{d^k}{dt^k}(e^{-x}x^{n+\alpha})$$

\noindent The generation function for generalized Laguerre polynomials is the following:

$$GF_{Laguerre}(t,w)=\sum_{k=0}^\infty L^{(\alpha)}_k(t)w^k = \frac{1}{(1-w)^{\alpha+1}}e^{-wt/(1-w)}$$

\noindent The generalized Laguerre polynomials are orthogonal over $[0,\infty]$ with respect to the measure with weighting function $t^\alpha e^{-t}:$

$$\int_0^\infty x^\alpha e^{-x}L^{(\alpha)}_n(t)L^{(\alpha)}_m(t)dt=\frac{\Gamma(n+\alpha+1)}{n!}\delta_{nm},$$
therefore, the orthonormal version of the generalized Laguerre polynomials will take the form
\begin{equation}
\label{P-laguerre}
\hat{L}^{(\alpha)}_k(t)=\sqrt{\frac{k!}{\Gamma(k+\alpha+1)}}x^{\alpha/2}e^{-x/2}L^{(\alpha)}_k(t).
\end{equation}
\noindent The upper bound for the sub-Gaussian norm of generalized Laguerre polynomials will take the form

\begin{equation}
\label{tau-laguerre}
\tau_\varphi(\xi_k)\leq\tau_{Laguerre} (w,k)= \sqrt{\frac{\Gamma(k+\alpha+1)}{k!}}\tau w^k
\end{equation}

\noindent Let us estimate the value of $C_N$  in this case. Following the same procedure as provided above for Legendre polynomials, with inequality (\ref{tau-laguerre}) in mind,

$$\tau_\varphi(\Delta_N(t))\leq  \sum_{k=0}^\infty \tau_\varphi(\xi_k) a_k(t) - \sum_{k=0}^N \tau_\varphi(\xi_k) \hat{a}_k(t)  $$
$$=\sum_{k=0}^\infty \tau_\varphi(\xi_k) \int_0^\infty f(t,\lambda) \hat{L}^{(\alpha)}_k(\lambda) d\lambda - \sum_{k=0}^N \tau_\varphi(\xi_k) \hat{a}_k(t) $$
$$= \int_0^\infty
\left(
\sum_{k=0}^{\infty}
\tau_\varphi(\xi_k)f(t,\lambda) \hat{L}^{(\alpha)}_k(\lambda) \right)d\lambda- \sum_{k=0}^N \tau_\varphi(\xi_k) \hat{a}_k(t) $$
$$\leq \left( \int_0^\infty |f(t,\lambda)|^2d\lambda\right)^{1/2}\left(\int_0^\infty \left(
\sum_{k=0}^{\infty}
\tau_\varphi(\xi_k) \hat{L}^{(\alpha)}_k(\lambda)\right)^2 d\lambda\right)^{1/2}- \sum_{k=0}^N \tau_\varphi(\xi_k) \hat{a}_k(t) $$
$$\leq  \left( \int_0^\infty|f(t,\lambda)|^2d\lambda\right)^{1/2}
$$
$$\times
\left(\int_0^\infty \left(
\sum_{k=0}^{\infty}
\sqrt{\frac{\Gamma(k+\alpha+1)}{k!}}\tau w^k \sqrt{\frac{k!}{\Gamma(k+\alpha+1)}}\lambda^{\alpha/2}e^{-\lambda/2}L^{(\alpha)}_k(\lambda)\right)^2 d\lambda\right)^{1/2}$$
$$- \sum_{k=0}^N \sqrt{\frac{\Gamma(k+\alpha+1)}{k!}}\tau w^k \hat{a}_k(t)  $$
$$= \tau \left( \int_0^\infty |f(t,\lambda)|^2d\lambda\right)^{1/2}\left(\int_0^\infty \lambda^{\alpha}e^{-\lambda}\left(
\sum_{k=0}^{\infty}
 w^k L^{(\alpha)}_k(\lambda)\right)^2 d\lambda\right)^{1/2}
 $$
 $$
 - \sum_{k=0}^N \sqrt{\frac{\Gamma(k+\alpha+1)}{k!}}\tau w^k \hat{a}_k(t)
 $$
$$= \tau \left( \int_0^\infty |f(t,\lambda)|^2d\lambda\right)^{1/2}\left(\int_0^\infty \lambda^{\alpha}e^{-\lambda}\frac{1}{(1-w)^{2\alpha+2}}e^{-2w\lambda/(1-w)} d\lambda\right)^{1/2}
$$
$$- \sum_{k=0}^N \sqrt{\frac{\Gamma(k+\alpha+1)}{k!}}\tau w^k \hat{a}_k(t)
$$
$$= \tau \left( \int_0^\infty |f(t,\lambda)|^2d\lambda\right)^{1/2}\left(\left.-\lambda^{\alpha+1} \left(\frac{1-w}{\lambda(1+w)}\right)^{\alpha+1}\Gamma\left(\alpha+1,\frac{\lambda(1+w)}{1-w}\right)\right|_0^\infty\right)^{1/2}
 $$
 $$- \sum_{k=0}^N \sqrt{\frac{\Gamma(k+\alpha+1)}{k!}}\tau w^k \hat{a}_k(t),$$

\noindent where $\Gamma(s,x)$ is the upper incomplete gamma function. Applying the integration boundaries, we obtain

$$\tau \left( \int_0^\infty |f(t,\lambda)|^2d\lambda\right)^{1/2}\left(\left.-\lambda^{\alpha+1} \left(\frac{1-w}{\lambda(1+w)}\right)^{\alpha+1}\Gamma\left(\alpha+1,\frac{\lambda(1+w)}{1-w}\right)\right|_0^\infty\right)^{1/2} $$
 $$- \sum_{k=0}^N \sqrt{\frac{\Gamma(k+\alpha+1)}{k!}}\tau w^k \hat{a}_k(t)
 $$
$$= \tau \left( \int_0^\infty |f(t,\lambda)|^2d\lambda\right)^{1/2}\left(\frac{1-w}{1+w}\right)^{\alpha+1}\Gamma(\alpha+1) - \sum_{k=0}^N \sqrt{\frac{\Gamma(k+\alpha+1)}{k!}}\tau w^k \hat{a}_k(t), $$
and this expression is defined for $\alpha>-1$.

\noindent As a result, we have

$$C_N\leq
{\int_0^T}
 \left(\tau \biggl( \int_0^\infty |f(t,\lambda)|^2d\lambda\right)^{1/2}\left(\frac{1-w}{1+w}\right)^{\alpha+1}\Gamma(\alpha+1)$$
  $$- \sum_{k=0}^N \sqrt{\frac{\Gamma(k+\alpha+1)}{k!}}\tau w^k \hat{a}_k(t)\biggr)^p d\mu(t).$$

Let us denote the right side of the last inequality by
$C_{N,Laguerre}$.

It follows from  Theorem \ref{lp_l2} and  Theorem \ref{lp_g2}, that the following statements are true.

\begin{theorem}
\label{Laguerre1}
 Let a stochastic process $X=\{X(t),t\in[0,T]\}$ belong to the space ${Sub}_\varphi (\Omega)$ with the Orlicz $N$-function
  $$\varphi(t)=\frac{t^\gamma}{\gamma},\,\,1<\gamma\leq 2,$$
  and let the process $X(t)$ admits the orthogonal decomposition  (\ref{mainseries})
 based on Laguerre orthonormal polynomial families (\ref{P-laguerre}).
Assume that $C_{N,Laguerre}<\infty$,
$\tau_\varphi (a_k(u))\leq \tau_\varphi (\hat{a}_k(u))$, $\forall u\in [-1,1]$ and condition \eqref{tau-laguerre} holds true.

\noindent The model (\ref{model}) $X_N(t)=\sum_{k=0}^N \xi_k \hat{a}_k(t)$ approximates the stochastic process $X(t)$
with given reliability $1-\alpha$ and accuracy $\delta$ in the space $L_p[0,T]$, if
$$\left\{\begin{array}{c}
C_{N,Laguerre}\leq \delta/(\beta \ln \frac{2}{\alpha})^{p/\beta},\\
 C_{N,Laguerre}<\delta /p^{p\left(1-1/\gamma\right)},
\end{array}\right.$$
where $\beta$ if a number that fits the condition $\frac{1}{\beta}+\frac{1}{\gamma}=1$.
\end{theorem}

\begin{theorem}
\label{Laguerre2}
 Let a stochastic process $X=\{X(t),t\in[0,T]\}$ belong to the space ${Sub}_\varphi (\Omega)$ with the Orlicz $N$-function
   $$\varphi(t)=\left\{\begin{array}{c}\frac{t^2}{\gamma}, t<1,\\
    \frac{t^\gamma}{\gamma},t\geq1,
   \end{array}\right.,$$
 where $\gamma>2$,
  and let the process $X(t)$ admits the orthogonal decomposition  (\ref{mainseries})
 based on Laguerre orthonormal polynomial families (\ref{P-laguerre}).
Assume that $C_{N,Laguerre}<\infty$,
$\tau_\varphi (a_k(u))\leq \tau_\varphi (\hat{a}_k(u))$, $\forall u\in[0;\infty)$ and condition \eqref{tau-laguerre} holds true.
 The model (\ref{model}) $X_N(t)=\sum_{k=0}^N \xi_k \hat{a}_k(t)$ approximates the stochastic process $X(t)$
with given reliability $1-\alpha$ and accuracy $\delta$ in the space $L_p[0,T]$, if
$$\left\{\begin{array}{c}
C_{N,Laguerre}\leq \delta/(\beta \ln \frac{2}{\alpha})^{p/\beta},\\
 C_{N,Laguerre}<\delta /p^{p\left(1-1/\gamma\right)},
\end{array}\right.$$
where $\beta$ if a number that fits the condition $\frac{1}{\beta}+\frac{1}{\gamma}=1$.
\end{theorem}

\subsection{Gegenbauer polynomials}

The Gegenbauer polynomials $C^{(\alpha)}_k(t)$ are defined for $t\in[-1; 1]$ as solutions of the Gegenbauer differential equation:

$$(1-t^2)y''-(2\alpha+1)ty'+k(k+2\alpha)y=0.$$
When $\alpha=1/2$, this equation reduces to the Legendre equation. Therefore, for the mentioned $\alpha=1/2$, the Gegenbauer polynomials are reduced to the Legendre polynomials.

\noindent The Gegenbauer polynomials can be represented as Gaussian hypergeometric series when this series is finite:

$$C_k^{(\alpha)}(t) = \frac{(2\alpha)_k}{k!}{}_2F_1\left(-k,2\alpha+k;\alpha+\frac{1}{2};\frac{1-t}{2}\right),$$
where ${}_2F_1(a,b;c;z)$ is a special function which is a solution of the Euler hypergeometric differential equation:

$$t(1-t)\frac{d^2 y }{dt^2} + (c-(a+b+1)t)\frac{dy}{dt} - a b y =0.$$

\noindent This function can also be represented in the form of a series

$${}_2F_1(a,b;c;z) = \frac{\Gamma(c)}{\Gamma(a)\Gamma(b)}\sum_{n=0}^\infty \frac{\Gamma(a+n)\Gamma(b+n)}{\Gamma(c+n)}\frac{z^n}{n!}.$$

\noindent The Rodrigues formula for the Gegenbauer polynomials is the following:

$$C_k^{(\alpha)}(t)=\frac{(-1)^k}{2^k k!} \frac{\Gamma(\alpha+\frac{1}{2})\Gamma(k+2\alpha)}{\Gamma(2\alpha)\Gamma(\alpha+k+\frac{1}{2})}(1-t^2)^{1/2-\alpha} \frac{d^k}{dt^k}(1-t^2)^{k+\alpha-1/2}.$$

\noindent The generating function is given as

$$\sum_{k=0}^\infty C_k^{(\alpha)} (t)w^k = \frac{1}{(1-2tw+w^2)^\alpha}.$$

\noindent  The Gegenbauer polynomials are orthogonal:

$$\int_{-1}^1 (1-t^2)^{\alpha-1/2}C_n^{(\alpha)} (t)C_m^{(\alpha)} (t)dt=\frac{\pi 2^{1-2\alpha}\Gamma(k+2\alpha)}{k! (k+\alpha) \Gamma^2 (\alpha)}\delta_{nm}.$$

\noindent Therefore, the orthonormal version of the generalized Gegenbauer polynomials is
\begin{equation}
\label{P-gegenbauer}
\hat{C}^{(\alpha)}_k(t)=\frac{\Gamma (\alpha)\sqrt{k! (k+\alpha)}}{\sqrt{\pi 2^{1-2\alpha}}\sqrt{\Gamma(k+2\alpha)}}{\sqrt{(1-t^2)^{\alpha-1/2}}}C^{(\alpha)}_k(t).
\end{equation}
\noindent We can assume that
\begin{equation}
\label{tau-gegenbauer}
\tau_\varphi(\xi_k)\leq\tau_{Gegenbauer} (w,k)= \sqrt{\frac{k! (k+\alpha)}{\Gamma(k+2\alpha)}}\tau w^k
\end{equation}

\noindent Let us estimate the value of $C_N$. We have

$$\tau_\varphi(\Delta_N(t))\leq  \sum_{k=0}^\infty \tau_\varphi(\xi_k) a_k(t) - \sum_{k=0}^N \tau_\varphi(\xi_k) \hat{a}_k(t)
$$
$$=\sum_{k=0}^\infty \tau_\varphi(\xi_k) \int_{-1}^1 f(t,\lambda) \hat{C}^{(\alpha)}_k(\lambda) d\lambda - \sum_{k=0}^N \tau_\varphi(\xi_k) \hat{a}_k(t)
$$
$$= \int_{-1}^1 \left(
\sum_{k=0}^{\infty}
\tau_\varphi(\xi_k)f(t,\lambda) \hat{C}^{(\alpha)}_k(\lambda) \right)d\lambda- \sum_{k=0}^N \tau_\varphi(\xi_k) \hat{a}_k(t)
$$
$$\leq \left( \int_{-1}^1 |f(t,\lambda)|^2d\lambda\right)^{1/2}\left(\int_{-1}^1 \left(
\sum_{k=0}^{\infty}
\tau_\varphi(\xi_k) \hat{C}^{(\alpha)}_k(\lambda)\right)^2 d\lambda\right)^{1/2}- \sum_{k=0}^N \tau_\varphi(\xi_k) \hat{a}_k(t)
$$
$$\leq  \left( \int_{-1}^1 |f(t,\lambda)|^2d\lambda\right)^{1/2}\left(\int_{-1}^1 \left(
\sum_{k=0}^{\infty}
\frac{\Gamma (\alpha)}{\sqrt{\pi 2^{1-2\alpha}}}\tau w^k{\sqrt{(1-\lambda ^2)^{\alpha-1/2}}}C^{(\alpha)}_k(\lambda)\right)^2 d\lambda\right)^{1/2}
$$
$$- \sum_{k=0}^N \sqrt{\frac{k! (k+\alpha)}{\Gamma(k+2\alpha)}}\tau w^k \hat{a}_k(t)
$$
$$=   \frac{\tau\Gamma (\alpha)}{\sqrt{\pi 2^{1-2\alpha}}}
\left( \int_{-1}^1 |f(t,\lambda)|^2d\lambda\right)^{1/2}
\left(\int_{-1}^1 (1-\lambda ^2)^{\alpha-1/2}\left(\sum_{k=0}^N  w^k C^{(\alpha)}_k(\lambda)\right)^2 d\lambda\right)^{1/2}
$$
$$- \sum_{k=0}^N \sqrt{\frac{k! (k+\alpha)}{\Gamma(k+2\alpha)}}\tau w^k \hat{a}_k(t)
$$
$$=  \frac{\tau\Gamma (\alpha)}{\sqrt{\pi 2^{1-2\alpha}}}\left( \int_{-1}^1 |f(t,\lambda)|^2d\lambda\right)^{1/2}\left(\int_{-1}^1 \frac{(1-\lambda ^2)^{\alpha-1/2}}{(1-2\lambda w+w^2)^{2\alpha}} d\lambda\right)^{1/2}-$$ $$- \sum_{k=0}^N \sqrt{\frac{k! (k+\alpha)}{\Gamma(k+2\alpha)}}\tau w^k \hat{a}_k(t)
$$
$$=  \frac{\tau\Gamma (\alpha)}{\sqrt{\pi 2^{1-2\alpha}}}
\left( \frac{1}{(2 \alpha -1) w} 4^{-\alpha } \left(1-\lambda ^2\right)^{\alpha - 1/2}    \left(\frac{(\lambda -1) (\lambda +1) w^2}{(w+1)^2(w-1)^2}\right)^{1/2-\alpha }  \left(w^2-2    \lambda  w+1\right)^{1-2 \alpha }\right.$$
$$\times \left.\left. F_1\left(1-2 \alpha ;\frac{1}{2}-\alpha    ,\frac{1}{2}-\alpha ;2-2 \alpha ;\frac{w^2-2 \lambda  w+1}{(w+1)^2},\frac{w^2-2    \lambda  w+1}{(w-1)^2}\right)  \right|_{-1}^1 \right)^{1/2}
   $$
   $$
   \times
   \left( \int_{-1}^1 |f(t,\lambda)|^2d\lambda\right)^{1/2}- \sum_{k=0}^N \sqrt{\frac{k! (k+\alpha)}{\Gamma(k+2\alpha)}}\tau w^k \hat{a}_k(t)
   $$
$$=  \frac{\tau\Gamma (\alpha)}{\sqrt{\pi 2^{1-2\alpha}}}
\left( \sqrt{\pi } \Gamma \left(\alpha +\frac{1}{2}\right) \left(w^2+1\right)^{-2
   \alpha } \, _2\tilde{F}_1\left(\alpha ,\alpha +\frac{1}{2};\alpha +1;\frac{4
   w^2}{\left(w^2+1\right)^2}\right) \right)^{1/2}
   $$
   $$
   \times
   \left( \int_{-1}^1 |f(t,\lambda)|^2d\lambda\right)^{1/2} - \sum_{k=0}^N \sqrt{\frac{k! (k+\alpha)}{\Gamma(k+2\alpha)}}\tau w^k \hat{a}_k(t)
   $$
$$=  \frac{\tau\Gamma (\alpha)}{\left(w^2+1\right)^{
   \alpha }}\sqrt{\frac{\Gamma \left(\alpha +\frac{1}{2}\right)}{\sqrt{\pi} 2^{1-2\alpha}}}    \, _2\tilde{F}_1^{1/2}\left(\alpha ,\alpha +\frac{1}{2};\alpha +1;\frac{4
   w^2}{\left(w^2+1\right)^2}\right)\left( \int_{-1}^1 |f(t,\lambda)|^2d\lambda\right)^{1/2} $$
   $$- \sum_{k=0}^N \sqrt{\frac{k! (k+\alpha)}{\Gamma(k+2\alpha)}}\tau w^k \hat{a}_k(t), $$
where $\, _2\tilde{F}_1(a,b;c;z)$ is the regularized hypergeometric function $\, _2F_1(a,b;c;z)$

As a result, we have

$$C_N\leq
{\int_0^T}
\left(\frac{\tau\Gamma (\alpha)}{\left(w^2+1\right)^{
   \alpha }}\sqrt{\frac{\Gamma \left(\alpha +\frac{1}{2}\right)}{\sqrt{\pi} 2^{1-2\alpha}}}    \, _2\tilde{F}_1^{1/2}\left(\alpha ,\alpha +\frac{1}{2};\alpha +1;\frac{4
   w^2}{\left(w^2+1\right)^2}\right)\right.
   $$
   $$\left.
   \times\left( \int_{-1}^1 |f(t,\lambda)|^2d\lambda\right)^{1/2}-
    \sum_{k=0}^N \sqrt{\frac{k! (k+\alpha)}{\Gamma(k+2\alpha)}}\tau w^k \hat{a}_k(t)\right)^p d\mu(t).$$

Let us denote the right side of the last inequality by
$C_{N,Gegenbauer}$.

It follows from  Theorem \ref{lp_l2} and  Theorem \ref{lp_g2}, that the following statements are true.

\begin{theorem}
\label{Gegenbauer1}
 Let a stochastic process $X=\{X(t),t\in [0,T]\}$ belong to the space ${Sub}_\varphi (\Omega)$ with the Orlicz $N$-function
  $$\varphi(t)=\frac{t^\gamma}{\gamma},\,\,1<\gamma\leq 2,$$
  and let the process $X(t)$ admits the orthogonal decomposition  (\ref{mainseries})
 based on Gegenbauer orthonormal polynomial families (\ref{P-gegenbauer}).
Assume that $C_{N,Gegenbauer}<\infty$,
$\tau_\varphi (a_k(u))\leq \tau_\varphi (\hat{a}_k(u))$, $\forall u\in [-1,1]$ and condition \eqref{tau-gegenbauer} holds true.
 The model (\ref{model}) $X_N(t)=\sum_{k=0}^N \xi_k \hat{a}_k(t)$ approximates the stochastic process $X(t)$
with given reliability $1-\alpha$ and accuracy $\delta$ in the space $L_p[0,T]$, if
$$\left\{\begin{array}{c}C_{N,Gegenbauer}\leq \delta/(\beta \ln \frac{2}{\alpha})^{p/\beta},\\
 C_{N,Gegenbauere}<\delta /p^{p\left(1-1/\gamma\right)},\end{array}\right.$$
where $\beta$ if a number that fits the condition $\frac{1}{\beta}+\frac{1}{\gamma}=1$.
\end{theorem}

\begin{theorem}
\label{Legendre2}
 Let a stochastic process $X=\{X(t),t\in [0,T]\}$ belong to the space ${Sub}_\varphi (\Omega)$ with the Orlicz $N$-function
   $$\varphi(t)=\left\{\begin{array}{c}\frac{t^2}{\gamma}, t<1,\\  \frac{t^\gamma}{\gamma},t\geq1,\end{array}\right.,$$
 where $\gamma>2$,
  and let the process $X(t)$ admits the orthogonal decomposition  (\ref{mainseries})
 based on Gegenbauer orthonormal polynomial families (\ref{P-gegenbauer}).
Assume that $C_{N,Gegenbauer}<\infty$,
$\tau_\varphi (a_k(u))\leq \tau_\varphi (\hat{a}_k(u))$, $\forall u\in [-1,1]$ and condition \eqref{tau-gegenbauer} holds true.
 The model (\ref{model}) $X_N(t)=\sum_{k=0}^N \xi_k \hat{a}_k(t)$ approximates the stochastic process $X(t)$
with given reliability $1-\alpha$ and accuracy $\delta$ in the space $L_p[0,T]$, if
$$\left\{\begin{array}{c}C_{N,Gegenbauer}\leq \delta/(\beta \ln \frac{2}{\alpha})^{p/\beta},\\ C_{N,Gegenbauer}<\delta /p^{p\left(1-1/\gamma\right)},\end{array}\right.$$
where $\beta$ if a number that fits the condition $\frac{1}{\beta}+\frac{1}{\gamma}=1$.
\end{theorem}

\section{Discussion and Conclusions}

The results obtained in this paper substantially expand the applicability of the method introduced in \cite{KMM2015,Chapter,Mokliachuk2012,Mokliachuk2014, Mokliachuk2018}  by enabling the use of three major families of classical orthonormal polynomials -- Legendre, generalized Laguerre, and Gegenbauer -- in the construction and analysis of $\varphi$-sub-Gaussian stochastic process models. Unlike Hermite and Chebyshev systems, these families typically lack closed-form normalized generating functions, which had limited their use in earlier frameworks. By introducing polynomial-specific upper bounds for the $\varphi$-sub-Gaussian norm $\tau_\varphi(\xi_k)$, we circumvent this difficulty and derive explicit reliability and accuracy estimates for models in $L_p([0,T])$.

Expanding the polynomial family available for the proposed method increases applicability and supports more scenarios  here this modeling approach may be introduced. The weight function associated with each polynomial family can be matched to the statistical properties of the process (e.g., its marginal distribution or covariance structure). This matching improves approximation quality and convergence rates.

Each polynomial family brings distinct structural advantages. Legendre polynomials, with uniform weight on a bounded interval, are especially suited for processes defined on finite domains with relatively homogeneous behavior. Generalized Laguerre polynomials, supported on $[0,\infty)$ and associated with exponential weights, naturally align with processes exhibiting decaying trajectories, such as those found in queueing systems, reliability engineering, or finance. Gegenbauer polynomials introduce a shape parameter, enabling flexible adaptation to
processes with different smoothness or symmetry properties. This tunability makes the Gegenbauer basis particularly powerful for modeling processes whose covariance structure changes across the domain.

The general technique used here - replacing exact normalized generating functions with polynomial-specific upper bounds - opens a path toward extending the reliability and accuracy theory to even broader classes of orthonormal systems. In particular, Jacobi polynomials present a natural next step. However, their two-parameter structure and more complex weight function introduce additional analytic challenges. Future work will focus on developing appropriate bounds that accommodate these challenges, enabling an analogous extension of the method to the full Jacobi family.

Overall, the framework developed in this paper provides practical, computable bounds that allow practitioners to determine the number of series terms required to achieve prescribed reliability and accuracy targets. This enhances the utility of orthonormal polynomial expansions in stochastic modeling, especially in applications where analytical tractability of coefficients is limited and numerical approximation is unavoidable. The next stage of research will involve numerical implementations and case studies, further validating the theoretical results and demonstrating the method's effectiveness in realistic modeling scenarios.

\section{Declarations}

This research was conducted independently by the author and outside the scope of employment. No employer data, confidential information, intellectual property, or internal materials were used. The views expressed are solely those of the author. The author declares no competing interests. No funding was received for this work.

\end{document}